\documentstyle[12pt]{amsart}

\makeatletter                                            
\renewcommand{\@begintheorem}[2]{                        
\rm \trivlist \item [\hskip \labelsep {\bf #2\ \ #1.}]   
                                }                        
\makeatother                                             

\newcommand{\ts}{\vspace{\baselineskip}\noindent{\bf Proof.}$\;\;$}
\newcommand{\ZZ}{{\bf Z}}
\newcommand{\QQ}{{\bf Q}}
\newcommand{\RR}{{\bf R}}

\newcommand{\CC}{{\bf C}}

\newcommand{\PP}{{\bf P}}

\newcommand{\half}{{1/2}}
\begin{document}

\title[Half twists of  Hodge structures]
{Half twists of  Hodge structures of CM-type}

\author{Bert van Geemen}
\address{Dipartimento di Matematica, Universit\`a di Pavia, 
Via Ferrata 1, I-27100 Pavia, Italia}
\email{geemen@@dragon.ian.pv.cnr.it}

\begin{abstract}
To a Hodge structure $V$ of weight $k$ with CM by a field $K$ we associate
Hodge structures $V_{-n/2}$ of weight $k+n$ for $n$ positive and, under 
certain circumstances, also for $n$ negative. We show that these `half twists' 
come up naturally in the Kuga-Satake varieties of weight two Hodge structures
with CM by an imaginary quadratic field. 
\end{abstract}

\maketitle
\begin{center}
{\sc Introduction}
\end{center}
\vskip15pt

A Hodge structure of CM-type is a Hodge structure $V$ on which a CM-field
$K$ acts. Given a CM-type (a set of certain complex embeddings of $K$),
we give a simple construction of Hodge structures $V_{-n/2}$ 
of weight $k+n$, for any positive integer $n$ where $k$ is the weight of $V$.
These half twists of $V$ are {\em not} related to the
Tate twists $V(-n)$ of $V$: $V_{-n}\not\cong V(-n)$. 

In certain circumstances one can also
define the half twist $V_\half$, of weight $k-1$. 
The geometry underlying the half twist in the case of 
hypersurfaces in projective space is investigated in  \cite{vGI}.
A basic case is Example \ref{cycn} in this paper.

In the first section we recall the basic definitions of Hodge structures
and we define the half twist.
In the second section we look at half twists from the point of view of representations of $\CC^*$. We also give a geometrical interpretation of
the half twist in case $V$ is a sub-Hodge structure of $H^k(X,\QQ)$
for a smooth projective variety $X$.
In the last section we consider Hodge structures $V$ of weight two with $\dim V^{2,0}=1$ and with a suitable action of an imaginary quadratic field. 
We show that the first cohomology group of the
Kuga-Satake variety associated to such a Hodge structure
has a summand which is the half twist of $V$ and that half twists can be used to understand the other summands as well. So, in a certain sense,
half twists `partly generalize' the Kuga-Satake construction which associates
weight one Hodge structures to certain weight two Hodge structures.
These results were motivated by an example of C.\ Voisin \cite{Vks} which also
inspired the general definition of half twist.

I'm indebted to E.\ Izadi and M.\ Nori for simplifications and generalizations of my original construction.

\section{Definitions and basic properties}

\subsection{Hodge structures.} 
Recall that a (rational) Hodge structure of weight $k\;(\in\ZZ_{\geq 0})$
is a $\QQ$-vector space $V$ with a decomposition of its complexification $V_\CC:=V\otimes_\QQ\CC$ (where complex conjugation is given by $\overline{v\otimes z}:=v\otimes\bar{z}$ for $v\in V$ and $z\in\CC$):
$$
V_\CC=\oplus_{p+q=k}V^{p,q},\qquad{\rm and}\quad \overline{V^{p,q}}=V^{q,p},
\qquad(p,\,q\in\ZZ_{\geq 0}).
$$
Note that we insist on $p$ and $q$ being non-negative integers throughout this paper, so we only consider `effective' Hodge structures.

\subsection{Hodge structures of CM-type.}\label{hodaut}
Let $V$ be a Hodge structure such that $V$ is also a $K$-vector space
for some $CM$-field $K$ and such that the Hodge decomposition on $V$ is stable
under the action of $K$:
$$
xV^{p,q}\,\subset\,V^{p,q},\qquad (x\in K,\;p,\,q\in\ZZ_{\geq 0}).
$$
In particular, $K\hookrightarrow End_{Hod}(V)$. We will then say that $V$ is a Hodge structure of CM-type (with field $K$).

\subsection{CM-types.}
To define the half twist we need to fix a CM-type of the field $K$.
Recall that the CM-field $K$ has $2r=[K:\QQ]$ complex embeddings $K\hookrightarrow \CC$ and that a CM-type is a subset $\Sigma=\{\sigma_1,\ldots,\sigma_r\}$ of distinct embeddings with the property that no two are complex conjugate. Hence if we define as us
ual
$\overline{\Sigma}:=\{\bar{\sigma}_1,\ldots,\bar{\sigma}_{r}\}$ then
any embedding of $K$ is either in $\Sigma$ or in $\overline{\Sigma}$.

\subsection{The half twists.}
Let $V$ be a Hodge structure of CM-type with field $K$.
We consider the eigenspaces of the $K$-action on the $V^{p,q}$'s:
$$
V^{p,q}_\sigma:=\{v\in V^{p,q}:\;xv=\sigma(x)v\quad \forall x\in K\},
\qquad \sigma:K\hookrightarrow \CC.
$$
Given a CM-type $\Sigma$, we define two subspaces of 
$V^{p,q}$ whose direct sum is $V^{p,q}$:
$$
V^{p,q}_\Sigma:=\oplus_{\sigma\in\Sigma} V^{p,q}_\sigma,\qquad 
V^{p,q}_{\overline{\Sigma}}:=\oplus_{\sigma\in\overline{\Sigma}} V^{p,q}_\sigma.
$$

The negative half twist of $V$, denoted by $V_{-\half}\;(=V_{\Sigma,-\half})$, 
is the decomposition 
of $V_\CC$ given by the subspaces:
$$
V_{-\half}^{r,s}:=V^{r-1,s}_{\Sigma}\oplus V^{r,s-1}_{\overline{\Sigma}}.
$$
It is not hard to see that this is a Hodge structure, of CM-type with field $K$, on $V$ of weight $k+1$ where $k$ is the weight of $V$. 
By successively performing
the negative half twist one obtains $V_{-n/2}$, a Hodge structure on $V$ of weight $k+n$.

We observe that to define the (positive) half twist one would put:
$$
V_\half^{r,s}:=V^{r+1,s}_\Sigma\oplus V^{r,s+1}_{\overline{\Sigma}},
$$
however, now the subspaces $V^{k,0}_{\overline{\Sigma}}$ and $V^{0,k}_{\Sigma}$ of $V_\CC$
do not appear in $(V_\half)_\CC$ and therefore this definition does not define a
Hodge structure on $V$ (and in general not on any $\QQ$-subspace of $V$).
In particular, we will define the half twist of $V$ only if $V^{k,0}_{\overline{\Sigma}}=0$
(the complex conjugate of this space is $V^{0,k}_\Sigma$ which is then also trivial).

\section{Half twists via representations}

\subsection{} Hodge structures can also be defined via representations
of $\CC^*$ (more precisely, algebraic representations of $Res_{\CC/\RR}(G_m)$
in $GL(V_\RR)$).
We determine the representations corresponding to the half twists. Proposition
\ref{halftens} often points out interesting geometry, as in example \ref{cycn}.
In \cite{vGI} more such examples are investigated.
 
\subsection{}
The Hodge structure on a $\QQ$-vector space $V$ defined by an
algebraic representation $h:\CC^*\rightarrow GL(V_\RR)$ will be denoted by $(V,h)$,
its Hodge decomposition is:
$$
V^{p,q}:=\{v\in V_\CC:\;h(z)v=z^p\bar{z}^qv\,\}.
$$
The usual algebra constructions on representations can be applied to Hodge structures.
In particular, given rational Hodge structures $(V,h)$, $(W,h_W)$ of weight $k,\,k_W$
their tensor product is 
the rational Hodge structure 
of weight $k+k_W$ defined by:
$$
h\otimes h_W:\CC^*\longrightarrow GL(V_\RR\otimes W_\RR),\qquad
z\longmapsto[v\otimes w\mapsto (h(z)w)\otimes (h_W(z)w)].
$$

\subsection{} Let $(V,h)$ be a Hodge structure of weight $k$ of CM-type with field $K$.
To identify the $\CC^*$-representation on $V_\RR$
which defines the half twist, consider the $\RR$-linear extension of
the action of $K$ on $V_\RR$. This gives an action of
$K\otimes_\QQ\RR$ on $V_\RR$ and recall that
$K\otimes_\QQ\RR\cong \oplus_{i=1}^r \CC$ (if you write 
$$
K=\QQ[X]/(f)\qquad{\rm then}\qquad K\otimes\RR\cong \RR[X]/(f)\cong
\prod_{i=1}^r \RR[X]/(f_i),
$$
where $f=\prod f_i$ is the decomposition of $f$ in irreducible polynomials
$f_i$ of degree 2 in $\RR[X]$).
Choosing a CM-type
$\Sigma=\{\sigma_1,\ldots,\sigma_r\}$ of $K$ specifies an isomorphism 
$$
K\otimes\RR\stackrel{\cong}{\longrightarrow} \oplus_{j=1}^r \CC,\qquad x\otimes t\longmapsto
(t\sigma_1(x),\ldots,t\sigma_r(x)),
$$ 
we will denote the $\RR$-linear extensions of
the $\sigma_i$'s by the same symbols. 
The inverse of this isomorphism is denoted by $\phi$
$$
\phi:\oplus_{j=1}^r \CC\stackrel{\cong}\longrightarrow K\otimes\RR,
\qquad{\rm with}\quad
\sigma_i(\phi(z_1,\ldots,z_r))=z_i\quad(1\leq i\leq r).
$$
Denoting by $x\mapsto\bar{x}$ the complex conjugation on $K$ as well as its $\RR$-linear extension to $K\otimes_\QQ\RR$ let $\sigma_{r+i}(x):=\overline{\sigma_i(x)}=\sigma_i(\bar{x})$, hence
$$
\sigma_{r+i}(\phi(z_1,\ldots,z_r))=\bar{z}_i\qquad(1\leq i\leq r).
$$

We define a homomorphism (depending on the CM-type $\Sigma$) by composing the
diagonal inclusion $\Delta$ with the isomorphism $\phi=\phi_\Sigma$:
$$
g_\Sigma:\CC^*\stackrel{\Delta}{\longrightarrow}\prod_{i=1}^r \CC^*
\stackrel{\phi}{\longrightarrow} (K\otimes_\QQ\RR)^*
$$
Since $(K\otimes_\QQ\RR)^*$ acts on $V_\RR$, we get a representation,
also denoted by $g_\Sigma$ of $\CC^*$ on $V_\RR$. As $V$ is of CM-type, the actions of $K$ and $h(\CC^*)$ commute and thus
the action of $g_\Sigma(\CC^*)\;(\subset (K\otimes_\QQ\RR)^*)$ on $V_\RR$ commutes with the action of $h(\CC^*)$.

To understand the action of $g_{\Sigma}$ we observe that
the eigenvalues of
$x\in K$ on $V_\CC$ are the $\sigma_i(x)$, $1\leq i\leq 2r$, each with the same multiplicity. The decomposition in eigenspaces is
$$
V_\CC=\oplus_{\sigma\in \Sigma\cup\overline{\Sigma}} V_{\CC,\sigma}\qquad
{\rm with}\quad xv=\sigma(x)v\qquad(x\in K\otimes_\QQ\RR,\;v\in V_{\CC,\sigma}).
$$
In particular, if $x=g_\Sigma(z)$ then since $g_\Sigma(z)=\phi(\Delta(z))=\phi(z,\ldots,z)$ 
we get:
$$
g_\Sigma(z)v=\phi(z,\ldots,z)v=\left\{\begin{array}{cl} 
zv\,&\;v\in \oplus_{\sigma\in\Sigma} V_{\CC,\sigma},\\
\bar{z}v&\;v\in \oplus_{\sigma\in\overline{\Sigma}} V_{\CC,\sigma},
\end{array}\right.
$$
This leads to an alternative, but equivalent, definition of the half twist:

\subsection{Definition.}
Let $(V,h,K)$ be a Hodge structure of CM-type and let $\Sigma$ be a CM-type of $K$.
For $n\in \ZZ$, the $n$-th half twist $V_{-n/2}$ is the $\CC^*$-representation defined by the homomorphism
$$
g_\Sigma^nh:\CC^*\longrightarrow GL(V_\RR),\qquad z\longmapsto g_\Sigma^{n}(z)h(z).
$$

\subsection{} For positive $n$ the representation space $V_{-n/2}$ is a Hodge structure of weight $k+n$ where $k$ is the weight
of $V$. Moreover, with the same CM-type,  one has:
$$
\left(V_{n/2}\right)_{m/2}=V_{(n+m)/2}.
$$
We already observed that the half twist $V_\half$ is a Hodge structure iff the eigenspaces for
the $K$-action
$V^{k,0}_\sigma$ are trivial for $\sigma\in \overline{\Sigma}$.

\subsection{Tate twists.}
The Tate Hodge structure $\QQ(n)$ ($n\in \ZZ$) is defined by
the vector space $\QQ$ and the homomorphism:
$$
h_n:\CC^*\longrightarrow GL_1(\RR),\qquad z\longmapsto (z\bar{z})^{-n},
$$
it has weight $-2n$ and $\QQ(n)^{p,q}=0$ unless $p=q=-n$
in which case $\QQ(n)^{-n,-n}=\CC$. It is convenient to allow negative weights for this Hodge structure.
The $n$-th Tate twist of $V$ is defined by $V(n):=V\otimes \QQ(n)$,
it is a Hodge structure of weight $k-2n$ with $V(n)^{p,q}=V^{p+n,q+n}$.

Since the homomorphism $h_n$ acts via scalar multiplication on $V_\RR$,
it commutes with the representation $g_\Sigma^mh$ which defines $V_{m/2}$.
Therefore one has
$$
\left(V_{m/2}\right)(n)=\left(V(n)\right)_{m/2}.
$$
This is also easy to see from the Hodge decompositions. Note that $V(-1)\neq 
(V_{-\half})_{-\half}$ since the representations $h_1$ and $g_\Sigma^2$ are
not equivalent.

\subsection{} \label{refhK}
The `abstractly' defined half twists are in fact sub-Hodge structures of rather natural Hodge structures of CM-type.
Recall that to the CM-field $K$ and the CM-type $\Sigma$
one can associate a weight one Hodge structure on the $\QQ$-vector space $K$. These Hodge structures, and the abelian varieties associated to them, have been extensively investigated,
cf. \cite{Lang}, \cite{DMOS}.
If we endow $K$ with the trivial Hodge structure of
weight zero, $K_\CC=K^{0,0}$, this Hodge structure is just $K_{-\half}$.
One can identify $K_{-\half}$ with $H^1(A_K,\QQ)$
for an abelian variety $A_K$ with CM by the field $K$ and whose CM-type is the one 
used to define the (negative) half twist.

Let $V$ be a Hodge structure of CM-type with field $K$ of weight $n$, then
$V\otimes_\QQ K_{-\half}$ is a Hodge structure of weight $n+1$ which has an action of the algebra $K\otimes_\QQ K$. The element $x\otimes y\in K\otimes K$ acts as $(x\otimes y)(v\otimes z)= xv\otimes yz$ on $V\otimes_\QQ K_{-\half}$.
We will identify some subspaces of $V\otimes_\QQ K_{-\half}$ with half twists of $V$.

\subsection{Proposition.}\label{halftens}
{\it
Let $(V,h,K)$ be a Hodge structure of CM-type. Then we have
an inclusion of Hodge structures:
$$
V_{-\half}\subset V\otimes_\QQ K_{-\half}, 
$$
more precisely:
$$
V_{-\half}=\left\{
w\in V\otimes_\QQ K_{-\half}:\; (x\otimes 1)w=(1\otimes x)w,\quad 
\forall x\in K
\right\}.
$$
If $V$ admits a half twist, $V_\half(-1)$ is also 
a sub-Hodge structure of $V\otimes_\QQ K_{-\half}$:
$$
V_\half(-1)=\left\{
w\in V\otimes_\QQ K_{-\half}:\; (x\otimes 1)w=(1\otimes \bar{x})w,\quad 
\forall x\in K
\right\}
$$
and the Hodge structure $V$ can be recovered from $V_\half$ by applying the first result to $V_\half$ rather then $V$:
$$
V\subset V_\half\otimes_\QQ K_{-\half}.
$$
}

\ts
We give two proofs.

The field $K$ acts on the $\QQ$-vector space $V$
and the (complex) eigenvalues of the action of $x\in K$ are the $\sigma_j(x)\in\CC$,
 each with the same multiplicity. In particular, $x,\,\bar{x}\in K$ are eigenvalues of the $K$-linear extension of the action of $x$ on $V\otimes_\QQ K$
and we denote by $V'$ and $V''$ the corresponding eigenspaces: 
$$
V\otimes_\QQ K\cong V'\oplus V''\oplus W,\qquad {\rm with}\quad
(x\otimes 1)w=\left\{\begin{array}{cc}\;(1\otimes x)w&\quad w\in V',\\
\;(1\otimes \bar{x})w&\quad w\in V'', \end{array}\right.
$$ 
and $W$ is a $K$-stable complementary subspace.
The projections on the summands give isomorphisms of $K=K\otimes 1$-vector spaces $V=V\otimes 1\rightarrow V_1$ and $V\rightarrow V_{r+1}$. 
These $V_j$ are sub-Hodge structures of $V\otimes K_{-\half}$ since the actions
of $K$ and $h$, $g$ commute. The Hodge structure on $K_{-\half}$ is defined by the $g(z)\in (K\otimes \RR)^*$ and thus $1\otimes g(z)=g(z)\otimes 1$ 
on $V'$ and $1\otimes g(z)=g(\bar{z})\otimes 1$ on $V''$.
Note also that $g(tz)=tg(z)$ for $t\in\RR$.
Therefore
$$
h(z)\otimes g(z)=h(z)g(\overline{z})\otimes 1=z\bar{z}h(z)g(z^{-1})\otimes 1
$$
on $V''$ which identifies the Hodge structure on $V''$ with the one
on $V_\half(-1)$.
The proof of $V'\cong V_{-\half}$ is similar but easier. 
The last result follows by applying the first result to $V_\half$ rather than $V$ and using $(V_\half)_{-\half}=V$, note this identifies $V$ with a specific subspace of $V_{\half}\otimes K_{-\half}$.

Another proof of the first result is as follows: since for $x\in K$ the map
$w\mapsto (1\otimes x-x\otimes 1)w$ is a $\QQ$-linear endomorphism of $V\otimes_\QQ K_{-\half}$, its kernel is a $\QQ$-vector space and therefore
$$
V':=\left\{
w\in V\otimes_\QQ K_{-\half}:\; (x\otimes 1)w=(1\otimes x)w,\quad 
\forall x\in K
\right\}
$$
is a $\QQ$-vector space. Its complexification $V'_\CC$ 
is the folllowing subspace of $(V\otimes K_{-\half})_\CC$:
$$
V'_\CC=\oplus_{\sigma\in \Sigma\cup\bar{\Sigma}}\, V_\sigma\otimes K_\sigma
$$
with $V_\sigma$ (resp.\ $K_\sigma$) the subspace of $V_\CC$ (resp.\ $K_\CC$) on which $x\in K$ acts as $\sigma(x)\;(\in\CC)$.
Since $K_{-\half}^{1,0}=\oplus_{\sigma\in\Sigma}K_\sigma$ we get:
$$
\left(V'\right)^{r,s}=(\oplus_{\sigma\in\Sigma} V_\sigma^{r-1,s}\otimes K_\sigma) \oplus (\oplus_{\sigma\in\bar{\Sigma}}
V_{{\sigma}}^{r,s-1}\otimes K_{{\sigma}})
$$
which, since $K_\sigma\cong\CC$, is just the definition of $V^{r,s}_{-\half}$
hence $V'\cong V_{-\half}$. The other statements can be proved in a similar fashion.
\qed

\subsection{Geometrical version of the half twist.}
The proposition shows that if $V$ is a sub-Hodge structure of $H^k(X,\QQ)$
for some projective variety $X$, then $V_\half(-1)$ is a sub-Hodge structure
of $H^k(X,\QQ)\otimes H^1(A_K,\QQ)$, which is itself a summand of 
$H^{k+1}(X\times A_K,\QQ)$. 

\subsection{Polarizations.} 
A polarization on the Hodge structure $(V,h)$ of weight $k$ is a bilinear map:
\[ 
\Psi: V \times V \longrightarrow \QQ  
\]
satisfying (for all $v,w\in V_{\RR}$): 
\[
\Psi(h(z)v,h(z)w)=(z\bar{z})^k\Psi(v,w)
\]
and 
\[ 
\Psi (v,h(i)w)\;\;\;
\mbox{\rm is a symmetric and positive definite form:}
\] 
$\Psi(v,h(i)w)=\Psi(w,h(i)v)$ for all $v,w\in V_\RR$ and 
$\Psi(v,h(i)v)>0$ for all $v\in V_\RR-\{0\}$.
The primitive cohomology groups of algebraic varieties are polarized.

A polarized Hodge structure of CM-type with field $K$ is a polarized Hodge structure $(V,h,\psi)$ such that $(V,h,K)$ is of CM-type and such that
$$
\Psi(xv,w)=\Psi(v,\bar{x}w),\qquad x\in K,\;v,\,w\in V.
$$

\subsection{A polarization on the half twist.}
The Hodge structure $K_{-\half}$ has a polarization, cf.\ \cite{Lang}, \cite{DMOS}.
If $V$ has a polarization, then also $V\otimes K_{-\half}$ has a polarization (the tensor product polarization) and by restriction one obtains a polarization on $V_\half$. 

One can also proceed more explicitly: if $\Psi$ is a polarization on
$V$ then one chooses an element $\alpha\in K$ such that 
$\bar{\alpha}=-\alpha$ and such that for $\sigma\in\Sigma$
the purely imaginary complex numbers
$\sigma(\alpha)$ all have positive imaginary part.
Then the bilinear form
$$
\Psi':V\times V\longrightarrow \QQ,\qquad\Psi'(v,w):=\Psi(v,\alpha w)
$$
is a polarization on $V_\half$. To verify all the properties it is
convenient to split 
$$
(V_\half)_\RR=\oplus_{i=1}^r (V_\half)_i
$$
corresponding to the decomposition $K\otimes_\QQ\RR=\oplus_{i=1}^r\CC$.

\subsection{Example.}\label{cycn}
Let $F\in\CC[X_0,\ldots,X_{n}]$, homogeneous of degree $n+2$, define a
smooth variety $Y\subset\PP^n$. Let 
$$
X=Zeroes(X_{n+1}^{n+2}-F)\qquad (\subset \PP^{n+1}),\qquad
{\rm let}\qquad V:=H^{n}(X,\QQ)_0
$$
be the primitive cohomology group of $X$.
Then $V$ is a vector space over the field $K$ of $n+2$-th roots of unity, 
where the roots of unity act by multiplication on the variable $X_{n+1}$.  

The vector space $H^{n,0}(X)$ is one dimensional, hence we can apply the half twist to $V$ (for any CM -type which includes the embedding
$\sigma$ of $K$ defined by $xv=\sigma(x)v$ for $v\in H^{n,0}(X)$).
Proposition \ref{halftens} implies that
$$
V=H^{n}(X,\QQ)_0\hookrightarrow V_\half\otimes K_{-\half},
$$
with $V_\half$ and $K_{-\half}$ Hodge structures of weight $n-1$ and $1$ respectively.
In this case it is not hard to see geometrically that such Hodge structures
exist.

Note that a general line $l\;\subset \PP^{n}$ meets the hypersurface
$Y:=Zeroes(F)$ (the branch locus of the natural map $X\rightarrow \PP^{n}$)
in $n+2$ distinct points, so these have $n+2-3=n-1$ moduli. Since the grassmanian of lines in $\PP^{n}$ has dimension $2(n-1)$, we find that for
a general set of $n+2$ points on $\PP^1$ there is a
$n-1$-dimensional family $S$ of lines $l\subset \PP^{n}$ each of which meet $Y$ in $n+2$ points with the same moduli. The union of these lines will be
$\PP^{n}$. Let $C_n$ be the inverse image of
(any) of these lines in $X$. After replacing $S$ by the desingularization
of a finite cover
$\tilde{S}$, we then get a dominant, rational map:
$$
\pi:\tilde{S}\times C_n\longrightarrow X.
$$
Since a rational map is defined outside a subset of codimension at least 2,
the pull-back of the regular $n$-form on $X$ extends to a regular
$n$-form on $\tilde{S}\times C_n$.
Hence we get $H^{n}(X,\QQ)_0\hookrightarrow H^{n-1}(\bar{S},\QQ)\otimes
H^1(C_n,\QQ)$ where $\bar{S}$ is some compactified desingularization of 
$\tilde{S}$. We will relate this example to Shioda's results on Fermat type hypersurfaces in \cite{vGI}.

\section{Kuga-Satake varieties}

\subsection{} \label{intro3}
To a polarized Hodge structure $(V,\psi)$ of weight 2 with 
$\dim V^{2,0}=1$ the construction
of Kuga and Satake associates a Hodge structure $(C^+(V),h_s)$
of weight 1 on the even Clifford algebra $C^+(V)$ of the quadratic space $(V,\psi)$ (see \cite{KS} and \cite{vG} for a detailed construction). 
It has the property that there is an inclusion of Hodge structures
$V\hookrightarrow C^+(V)\otimes C^+(V)$.
The (isogeny class of) abelian variety associated to $C^+(V)$ is called the Kuga-Satake variety
$KS(V)$ of $V$, so $H^1(KS(V),\QQ)=C^+(V)$.

In the remainder of this paper we consider such Hodge structures which are of
CM-type with an imaginary quadratic field $K$.
Since $\dim V^{2,0}=1$, the half twist $V_\half$ is a Hodge structure and has
weight one. 
Our main result is Theorem \ref{ksdec} which shows that
$V_\half$ is a summand of $(C^+(V),h_s)$. We also determine the other
summands and relate them to half twists.

This completes the results of C.\ Voisin in \cite{Vks}, she already found that 
two summands of $C^+(V)$, $S_0$ of dimension 2 and $S_1$ with $\dim S_1=\dim V$,
such that $V\hookrightarrow S_0\otimes S_1$. We will identify $S_1$ with $V_\half$ in Theorem \ref{ksdec}.
To find the simple summands of $(C^+(V),h_s)$ we use the Mumford-Tate group of
the Hodge structure $V$. 

\subsection{The Mumford-Tate group.}
Recall that the Special Mumford-Tate group $SMT(V)$ of a polarized Hodge structure 
$(V,h,\psi)$ is the smallest algebraic subgroup $G$ of $GL(V)$, defined over $\QQ$, for which
 $h(z)\in G(\RR)$ $(\subset GL(V)(\RR))$ for all $z\in \CC$ with $z\bar{z}=1$.

The simple summands of the Hodge structure $(C^+(V),h_s)$ are then the irreducible subrepresentations of $SMT(V)$ in $C^+(V)$.
It takes some rather long computations to determine these summands though.

The following lemma recalls the basic facts on the Mumford Tate group in 
this situation.

\subsection{Lemma.}\label{lem2}
{\it
Let $(V,h,\psi)$ be a weight 2 rational polarized Hodge structure of CM-type
by an imaginary quadratic field $K=\QQ(\phi)$, with $\bar{\phi}=-\phi$
and $\phi^2=-d$.

\begin{enumerate}
\item[i)]
The $\QQ$-bilinear map:
$$
H:V\times V \longrightarrow K,\qquad H(v,w):=\psi(v,w)+\phi^{-1}\psi(v,\phi w)
$$
is a hermitian form on the $K:=\QQ(\phi)$-vector space $V$ (so
$\overline{H(v,w)}=H(w,v)$ and $H$ is $K$-linear in the second variable).

\item[ii)]
We have
$$
SMT(V)\,\subset\, U(H):=\{g\in GL(V):\;H(gv,gw)=H(v,w),\qquad
g\phi=\phi g\},
$$
the unitary group of the $K$-vector space $V$ with the hermitian form $H$. 

\item[iii)]
In case $(V,h)$ has weight two, (the $\RR$-linear extension of)
the hermitian form $H$ has signature $(\half\dim V^{1,1},\dim V^{2,0})$ on the
complex vector space $V\otimes_\QQ\RR$.
\end{enumerate}
}

\ts
Since the weight is even, $\psi$ is symmetric and: 
$$
\begin{array}{rcl}
H(w,v)&=&\psi(w,v)+\phi^{-1}\psi(w,\phi v)\\
&=& \psi(v,w)+\phi^{-1}\psi(\phi v,w)\\
&=& \psi(v,w)+(d\phi)^{-1}\psi(\phi^2 v,\phi w)\\
&=& \psi(v,w)-\phi^{-1}\psi(v,\phi w)\\
&=& \overline{H(v,w)}.
\end{array}
$$
The $K$-linearity in the second factor follows from:
$$
\begin{array}{rcl}
H(v,\phi w)&=& \psi(v,\phi w)+\phi^{-1}\psi(v,\phi^2 w)\\
&=&\psi(v,\phi w)-d\phi^{-1}\psi(v, w)\\
&=&\phi\left( \psi(v,w)+\phi^{-1}\psi (v,\phi w)\right)\\
&=& \phi H(v,w).
\end{array}
$$
Using that 
$\psi(h(z)v,h(z)w)=(z\bar{z})^k\psi(v,w)$ where $k$ is the weight of $V$ and
that $\phi h(z)=h(z)\phi$ we get:
$$
\begin{array}{rcl}
H(h(z)v,h(z)w)&=& \psi(h(z)v,h(z)w)+\phi^{-1}\psi(h(z)v,\phi h(z)w)\\
&=& (z\bar{z})^k\psi(v,w)+\phi^{-1}\psi(h(z)v,h(z)\phi w)\\
&=& (z\bar{z})^k\left( \psi(v,w)+\phi^{-1}\psi(v,\phi w)\right)\\
&=&(z\bar{z})^k H(v,w),
\end{array}
$$
hence $SMT(V)\subset U(H)$.

In case the weight of $V$ is two, we write:
$$
V\otimes_\QQ\RR=V_1\oplus V_2\qquad{\rm with}\quad
V_1\otimes_\RR\CC=V^{1,1},\quad V_2\otimes_\RR\CC=V^{2,0}\oplus V^{0,2}.
$$
Then $\psi$ is negative definite on $V_2$ and positive definite on $V_1$.
Since $\phi$ and $h(z)$ commute, $\phi$ maps the eigenspaces $V^{p,q}$ into themselves hence $\phi(V_i)\subset V_i$. As $\phi^2=-d$ with $d>0$, $\phi$ does not have real eigenvalues and we can choose a $\RR$-basis $f_1,\phi f_1$,
$\ldots$,$f_r,\phi f_r$ of $V_2$ (and similarly for $V_1$).
Since $\psi(f_i,\phi f_i)=$ $d^{-1}\psi(\phi f_i,\phi^2 f_i)$
$=-\psi(\phi f_i,f_i)=-\psi(f_i,\phi f_i)$,
we may assume this basis to be orthonormal. 

Since $f_1,\ldots, f_r$ is a $\CC=K\otimes_\QQ\RR$-basis of $V_2$ and
$H(f_i,f_j)=\psi(f_i,f_j)$, we see that $H$ is negative definite on $V_2$
(and $H$ is positive definite on $V_1$).
\qed

\subsection{Remark.}\label{moduli}
If $(V,h,\psi)$ and $K=\QQ(\phi)$ are as in section\ref{intro3},
then the Lie group $U(H)(\RR)$ of real points of $U(H)$ is ismomorphic to
$U(1,m-1)$.

For any $g\in U(H)(\RR)$ the Hodge structure $(V,h^g)$ with 
$$
h^g:\CC^*\longrightarrow GL(V)(\RR),\qquad z\longmapsto gh(z)g^{-1}
$$
is also polarized by the same $\psi$ and has $\phi\in End_{Hod}(V,h^g)$.
This implies that $SMT(V,h^g)=U(H)$ for any general $g\in U(H)(\RR)$.
One can show that the moduli space of such Hodge structures is
isomorphic to the complex $m-1$-ball
$U(1,m-1)/(U(1)\times U(m-1))$.

The inclusion $U(H)\subset SO(\psi)$ induces $U(1,m-1)\subset SO(2,2m-2)$,
this well-known inclusion is used for example to restrict modular forms from othogonal groups to unitary groups.

\subsection{} The universal cover of the orthogonal group $SO(\psi)$
has a natural spin representation on the even Clifford algebra $C^+(V)$.
The following proposition describes the spin representation over $\QQ$
for the quadratic forms under consideration. Over the complex numbers these results are very well known, but over a number field the situation is a bit
delicate. The proposition will be  used to decompose the
spin representation as a representation of $U(H)$, the Mumford Tate group of a general Hodge structure of CM-type $(V,h,\psi,K)$ under consideration.

\subsection{Proposition.}\label{sodec}
{\it
 Let $(V,h,\psi)$ be a polarized weight $2$ Hodge structure with $\dim V^{2,0}=1$ which is of CM-type for an imaginary quadratic field $K$.
Then there is a $\QQ$-basis of $V$ such that
$$
\psi(x,x)=\sum_{i=1}^m d_ix_i^2+d\sum_{i=1}^m d_ix_{m+i}^2,\qquad 
{\rm with}\quad d_1<0,\quad d_2,\ldots,d_m>0.
$$

The representation of $so(2m)$ on $C^+(Q)$ decomposes as 
$C^+(V)=S^{2^{m-2}}$ where $S$ is a $so(\psi)$-representation of dimension $2^{m+1}$ whose irreducible components are:
$$
S\cong\left\{ \begin{array}{llll}
S_+\oplus S_-,\quad& End_{so(\psi}(S_\pm)\cong D,\quad&S_\pm\otimes\CC\cong\Gamma_\pm^2\qquad & {\rm if}\;m\equiv 0\,(4),\\
S_1^2,\quad& End_{so(\psi)}(S_1)\cong K,\quad&S_1\otimes\CC\cong\Gamma_+\oplus\Gamma_-& {\rm if}\;m\equiv 1,\,3\,(4),
\end{array}
\right.
$$
with $D$ a skew field of degree $4$ over $\QQ$,
and
$\Gamma_+$, $\Gamma_-$ are the two half-spin representations
of $so(2m)\otimes_\QQ\CC$, each of which has dimension $2^{m-1}$.

In case $m\equiv 2 (4)$ there are two possibilities. If the equation
$-\prod d_i=x^2+dy^2$ has a solution $(x,y)\in\QQ^2$, then 
$$
S\cong S_+^2\oplus S_-^2,\qquad End_{so(\psi)}(S_\pm)\cong \QQ,\qquad S_\pm\otimes\CC\cong\Gamma_\pm
$$ 
(the split case). In case this equation has no solution
we have 
$$
S\cong S_+\oplus S_-,\qquad End_{so(\psi)}(S_\pm)\cong D,\qquad S_\pm\otimes\CC\cong\Gamma_\pm^2
$$ 
(and $D$ is a skew field of degree $4$ over $\QQ$, we call this the non-split case).
}

\ts
To find this basis of $V$, choose $e_1\in V$ with $\psi(e_1,e_1)\neq 0$ and let
$d_1:=\psi(e_1,e_1)$,
$e_{m+1}:=\phi e_1$. Then $\psi(e_{m+1},e_{m+1})=d\psi(e_1,e_1)=dd_1$ and $d\psi(e_1,e_{m+1})=
\psi(\phi e_1,(-d)e_{1})=-d\psi(e_1,e_{m+1})$, hence $\psi(e_1,e_{m+1})=0$.
Next we take $e_2\in \langle e_1,\,e_{m+1}\rangle^\perp$ etc.
Since the signature of $\psi$ is $(2-,(2n-2)+)$, we may assume $d_1<0$ and
$d_2,\ldots,d_m>0$.

We recall that over $K$ the spin representation of $so(\psi)_K:=so(\psi)\otimes_\QQ K$ on $C^+(V)_K:=C^+(V)\otimes_\QQ K$
decomposes as a direct sum of $2^{m-1}$ copies of $\Gamma_+\oplus \Gamma_-$ and
then we take Galois invariants to find the irreducible $so(\psi)$-representations over $\QQ$.

Let $V_K:=V\otimes_\QQ \QQ(\sqrt{-d})$, 
and consider the following $K$-basis of $V_K$:
$$
f_i:=\mbox{$\frac{1}{2dd_i}$}(\sqrt{-d}e_i+e_{m+i}),\qquad 
f_{m+i}:=\half(-\sqrt{-d}e_i+e_{m+i}),\qquad(1\leq i\leq m)
$$
where we wrote $\sqrt{-d}e_i$ for $e_i\otimes\sqrt{-d}$. One verifies that:
$$
\psi(\sum y_jf_j,\sum y_jf_j)=\sum_{i=1}^m y_iy_{m+i}.
$$
Since $\psi(f_j,f_j)=0$, $\psi(f_j+f_k,f_j+f_k)=1$ if $|j-k|=m$ and is zero otherwise,  we have in $C(V)_K$:
$$
f_j^2=0,\quad f_if_{m+i}+f_{m+i}f_i=1,\qquad 
f_jf_k=-f_kf_j\quad{\rm if}\; |j-k|\neq 0,\,m.
$$
We denote conjugation on $K$ by a `$\bar{\phantom{x}}$', this acts on $C(V)_K$ via the second factor and:
$$
\bar{f}_i=\mbox{$\frac{1}{dd_i}$}f_{m+i},\qquad \bar{f}_{m+i}=d_idf_i.
$$
Let
$$
f:=f_{m+1}f_{m+2}\ldots f_{2m},\qquad{\rm then}\quad
 \bar{f}=\bar{f}_{m+1}\ldots\bar{f}_{2m}=(d^m\prod_id_i)f_1\ldots f_m.
$$
In $C(V)_K$ we have:
$$
f\bar{f}f=\delta f,\qquad \bar{f} f\bar{f}=\delta \bar{f},\qquad
{\rm with}\quad
\delta:=(-1)^{m(m-1)/2}d^m\prod_{i=1}^m d_i,
$$
the second is just the conjugate of the first which is an easy computation:
$$
\begin{array}{rcl}
f\bar{f}f&=& (d^m\prod_id_i)
(f_{m+1}\ldots f_{2m})(f_1\ldots f_m)(f_{m+1}\ldots f_{2m})\\
&=&
\delta (f_{m+1}\ldots f_{2m})(f_mf_{m-1}\ldots f_2f_1)(f_{m+1}\ldots f_{2m})\\
&=&
\delta
(f_{m+1}\ldots f_{2m})(f_m\ldots f_2)(1-f_{m+1}f_1)(f_{m+2}\ldots f_{2m})\\
&=&
\delta
(f_{m+1}\ldots f_{2m})(f_m\ldots f_3)(1-f_{m+2}f_2)(f_{m+3}\ldots f_{2m})-0\\
&=&\ldots\\
&=&\delta f.
\end{array}
$$
Now we consider the left $C(V)_K$-modules generated by $f$ and $\bar{f}$.
These modules are isomorphic, in fact the relations we just proved imply that
$$
R_{\bar{f}}:C(V)_Kf\longrightarrow C(V)_K\bar{f},\qquad xf\longmapsto xf\bar{f}
$$
is an isomorphism of left $C(V)_K$-modules with inverse
$$
R_f:C(V)_K\bar{f}\longrightarrow C(V)_K f,\qquad y\bar{f}\longmapsto y\bar{f}f.
$$
In \cite{FH}, Chapter 20, an inclusion $so(\psi)_K\hookrightarrow C^+(V)_K\,(\subset C(V)_K)$ is constructed and it is shown
(\cite{FH}, 20.19 and 20.20) that, as an $so(\psi)_K$-module, $C(V)_Kf$ is
isomorphic to the direct sum of the two half-spin representations:
$$
C(V)_Kf\cong \Gamma_{+,K}\oplus \Gamma_{-,K},\qquad
\Gamma_{\pm,K}\otimes_K\CC\cong \Gamma_\pm.
$$
Moreover, $C(V)_K^+$, the even Clifford algebra, is isomorphic to a product of
two matrix algebras:
$$
C(V)_K^+\cong M_{2^{m-1}}(K)\times M_{2^{m-1}}(K),
$$
this implies that the spin representation of $so(\psi)_K$ on $C^+(V)_K$
is isomorphic to $2^{m-1}$ copies of $\Gamma_{+,K}\oplus \Gamma_{-,K}$.

Before considering the situation over $\QQ$ we recall that the center of
$C^+(V)$ is $\QQ\oplus \QQ z$ with $z:=e_1e_2\ldots e_{2m}$ and that 
$z^2=(-1)^md^m\prod_i d_i^2$ (cf.\ \cite{Lam}, $\S$ 5.2). Since $\sqrt{(-1)^md^m\prod_i d_i^2}\in K$, the center of $C(V)^+_K$ is $K\times K$ and one can verify that
$\Gamma_{\pm,K}$ are the two eigenspaces of $z$ in $C(V)_Kf$.
We also observe that $C(V)_Kf\,\cap\, C(V)_K\bar{f}=\{0\}$, in fact if
$af=b\bar{f}$ then $af\bar{f}=b\bar{f}^2=0$, hence $0=af\bar{f}f=\delta af$,
hence $af=0$.

The subspace $C(V)_Kf$ is not defined over $\QQ$ in general, but the direct sum
$$
S_K:=\,C(V)_Kf\,\oplus\, C(V)_K\bar{f}
$$
obviously is defined over $\QQ$, that is $S_K=S\otimes_\QQ K$ for some $\QQ$-vector space $S\subset C(V)$, in fact
$$
S=C(V)(f+\bar{f})\,+\, C(V)\sqrt{-d}(f-\bar{f}).
$$
Moreover, $S$ is a representation space
for $so(\psi)\;(\subset so(\psi)_K)$.


To decompose $S$ into irreducible components
we determine $A:=End_{so(\psi)}(S)$, the endomorphisms of $S$ which commute with $so(\psi)$. Since $S_K\cong \Gamma_+^2\oplus\Gamma_-^2$, we have $A_K:=A\otimes_\QQ K\cong M_2(K)\times M_2(K)$, hence $\dim_\QQ A=8$.
It is clear that $A_K$ is generated by the center of $C^+(V)_K$ and the maps
$R_f$ and $R_{\bar{f}}$. To determine $A$ it suffices to find the invariants under conjugation in $A_K$.

Obviously the center of $C^+(V)$ lies in $A$. Moreover, the maps 
$$
\alpha,\;\beta:S\longrightarrow S,\qquad \alpha:x\mapsto x(f+\bar{f}),\qquad
\beta:x\longmapsto x\sqrt{-d}(f-\bar{f})
$$
commute with $so(\psi)$ (which acts from the left whereas $\alpha$ and $\beta$ act from the right). Note we have:
$$
(f+\bar{f})^2=f^2+f\bar{f}+\bar{f}f+\bar{f}^2=f\bar{f}+\bar{f}f=\delta,
$$
the last equality holds since in $S_K$ we have:
$$
(af+b\bar{f})(f\bar{f}+\bar{f}f)=0+af\bar{f}f+b\bar{f}f\bar{f}+0=
\delta(af+b\bar{f})
$$
and similarly $(\sqrt{-d}(f-\bar{f}))^2=d(f\bar{f}+\bar{f}f)=d\delta$.
Moreover, $\alpha\beta=-\beta\alpha$. Thus the $\QQ$-algebra generated by $\alpha$ and $\beta$ is the quaternion algebra $D:=(\delta,d\delta)$ and 
$$
A\cong D\otimes_\QQ\QQ(z),\qquad{\rm with}\quad D=(\delta,d\delta).
$$  
Since $(\alpha\beta)^2=-\alpha^2\beta^2=-d\delta^2$, $D$ contains a 
copy of the field $\QQ(\sqrt{-d})\cong K$.

The center of $C^+(V)$ is $\QQ(z)$ with $z^2=(-d)^m\prod_i d_i^2$, hence:
$$
\QQ(z)\cong \QQ\times \QQ\qquad{\rm if}\;m\equiv 0\;(2),\qquad
\QQ(z)\cong K\qquad{\rm if}\;m\equiv 1\;(2).
$$
As $d>0$, $d_1<0$ and
$d_2,\ldots,d_m>0$, the sign of $\delta$ is the sign of
$(-1)^{m(m-1)/2}(-1)$. Thus $\delta,\,d\delta$ are both negative if $m\equiv
0,\,1\;(4)$ so $D\otimes_\QQ\RR$ is isomorphic to the algebra of quaternions, a skew field, and thus $D$ is also a skew field. 

Hence in case $m\equiv 0\;(4)$ we have $A\cong D\times D$ with a skew field $D$,
therefore $S$ splits up in two components $S_+$ and $S_-$ with $End_{so(\psi)}(S_\pm)\cong D$ and since $D\otimes_\QQ\CC\cong M_2(\CC)$ the
$S_\pm\otimes_\QQ\CC$ are both direct sums of two copies of one irreducible $so(\psi)_\CC$ representation.

In case $m\equiv 1\;(4)$, $\QQ(z)\cong K$ and $D$ contains a copy of $K$ hence
$A\cong M_2(K)$ and thus $S$ is isomorphic to the sum of two isomorphic representations, irreducible over $\QQ$,
each of which, after tensoring with $K$, is the direct sum of two non-isomorphic representations.

In case $m\equiv 2\;(4)$, we have $A\cong D\times D$, but
the structure of $D$ depends on the $d_i$. In fact, $D\cong(-d,\delta)$ (using $\alpha\beta$ and $\alpha$ as generators) and
$\delta=-d_1d_2\ldots d_m d^{2k}$ (with $2k=m$), we also have 
$D\cong (-d,n)$ with $n:=-\prod_i d_i\in\ZZ_{>0}$. Hence $D\cong M_2(\QQ)$ iff
$n=x^2+dy^2$ for some $x,\,y\in\QQ^2$.   

In case $m\equiv 3\;(4)$ we have, as in the case $m\equiv 1\;(4)$, that
$A\cong M_2(K)$.
\qed

\subsection{}
For the general $(V,h,\psi,K)$ we consider, the Mumford Tate group is $U(H)$.
Thus the simple factors of the Hodge structure $(C^+(V),h_s)$ associated to $(V,h,\psi)$ are exactly the irreducible subrepresentations of the Lie algebra $u(H)$ ($\subset so(\psi)$) of $U(H)$ in $C^+(V)$ (which are defined over $\QQ$). 
We now determine the restriction of the $so(\psi)$ representation $S$ from Proposition \ref{sodec} to  $u(H)$.

\subsection{Proposition.}\label{sudec}
{\it
Let $S$ be the $2^{m+1}$-dimensional $so(\psi)$-representation defined in
Proposition \ref{sodec}. The $u(H)$-representation $S$ decomposes as follows:
$$
S\cong S_0\oplus S_1\oplus\ldots\oplus S_m,\qquad S_i\cong S_{m-i},\qquad
\dim_\QQ S_i=2{m\choose i}.
$$
The $S_i$ are irreducible $u(H)$-representations
except $S_l$ if $2l=m$, $l\equiv 2\;(4)$ and we are in the split case,
in that case $S_l\cong (S_l')^2$  and $S'_l$ is irreducible.

The $S_i$ are $K$-vector spaces and:
$$
End_{u(H)}(S_i)=\left\{\begin{array}{cc}K\;&\;(2i\neq m)\\ D\;&\;(2i=m),
\end{array}\right.
$$
with $D$ a quaternion algebra (a skew field except for the split case)
which contains $K$, but $End_{u(H)}(S'_l)=\QQ$.
}

\ts
First we determine the inclusion $u(H)_K\hookrightarrow so(\psi)_K$ and the 
restriction of the half-spin representations $\Gamma_{\pm,K}$ to $u(H)_K$.

Extending the scalars
from $\QQ$ to $K$, the endomorphism $\phi$ of $V$ has two eigenspaces in
$V\otimes_\QQ K$:
$$
V_K:=V\otimes_\QQ K=V_+\oplus V_-,\qquad (\phi\otimes 1)v=\pm(1\otimes \phi)v
\qquad (v\in V_\pm).
$$
Each of the eigenspaces is isotropic for the $K$-linear extension of $\psi$ to
$V_K$, in fact, 
$$
\begin{array}{rcl}
\psi(v,w)&=&d^{-1}\psi((\phi\otimes 1)v,(\phi\otimes 1)w)\\
&=&
d^{-1}\psi(\pm(1\otimes \phi)v,\pm(1\otimes \phi)w)\\
&=&d^{-1}\phi^2\psi (v,w)\\
&=&-\psi(v,w).
\end{array}
$$

The actions of $u(H)$ and $K$ on $V$ commute, hence $u(H)$ acts on the eigenspaces $V_+$ and $V_-$, each of which has dimension $m$. In particular we have a Lie algebra map
$u(H)\hookrightarrow gl(V_+)$, which, for dimension reasons, gives
an isomorphism: $u(H)_K\cong gl(V_+)$.
Since $\psi$ is preserved, this fixes the map $u(H)\rightarrow gl(V_-)$
(in fact $\psi$ gives a duality $V_-\stackrel{\cong}{\rightarrow} V_+^*$)
and we get (with respect to the basis $f_1,\ldots,f_{2m}$ of the proof of 
Prop.\ \ref{sodec}):
$$
u(H)_K\cong gl(V_+)\longrightarrow so(\psi)_K,\qquad
A\longmapsto \left(\begin{array}{cc} A&0\\ 0& -{}^tA\end{array}\right).
$$

This inclusion $gl(V_+)\hookrightarrow so(\psi)$ is the same as the one obtained from composing the isomorphism $gl(V_+)\cong end(V_+)\cong$ $V_+^*\otimes V_-\cong V_-\otimes V_+$ and the inclusion (cf.\ \cite{FH}, formula (20.4))
$V_-\otimes V_+\hookrightarrow \wedge^2 V_K\cong so(\psi)$.
>From \cite{FH}, 20.15 we obtain
the restrictions of the half-spin representations $\Gamma_{\pm,K}$ to $u(H)_K$:
$$
\Gamma_{+,K}\cong \oplus_{i} \Lambda^{2i} V_+,\qquad
\Gamma_{-,K}\cong \oplus_{i} \Lambda^{2i+1} V_+,
$$
(this is actually only an isomorphism of $su(H)_K$-representations).

To find the irreducible representations of $su(H)$
which are defined over $\QQ$ we need to know how the conjugation on $K$ acts. For this we use the following $C(V)$-modules (\cite{FH}, 20.12):
$$
C(V)f\cong \oplus_{i=0}^m \Lambda^i V_+,\qquad
C(V)\bar{f}\cong\oplus_{i=0}^m \Lambda^i V_-.
$$
The $f_{i_1}\ldots f_{i_r}f$ with
$1\leq i_1<\ldots <i_r\leq m$ are a basis of $C(V)_Kf$, 
and these correspond to the elements $f_{i_1}\wedge\ldots \wedge f_{i_r}
\in \Lambda^r V_+$. Their conjugate is, up to a constant, $f_{i_1+m}\ldots
f_{i_r+m}\bar{f}$ which corresponds to $f_{i_1+m}\wedge\ldots
\wedge f_{i_r+m}\in \Lambda^r V_-\cong \Lambda^{m-r} V_+$.
In particular, we get conjugation-invariant subspaces
$$
S_{i,K}= (\Lambda^i V_+)f\,\oplus\, (\Lambda^{m-i} V_+)\bar{f}\qquad
(\subset S_K=C(V)_Kf\oplus C(V)_K\bar{f})
$$
for $0\leq i\leq m$. Thus $S_{i,K}=S_i\otimes_\QQ K$ for a subspace $S_i\in C(V)$ which is invariant under $su(H)$.
If $i\neq m-i$, the two summands are not isomorphic as $gl(V_+)_K$-representations, hence $S_i$ is an irreducible $su(H)$-representation (if $W\subset S_i$ is an invariant subspace, then $W\otimes_\QQ K$ is a $gl(V_+)_K$ and Galois invariant subspace of $
S_{i,K}$
and hence $W=\{0\}$ or $W=S_i$). 
Since $End_{sl(V_+)}(S_{i,k})=K^2$ and $S_i$ is irreducible, $B:=End_{su(H)}(S_i)$ is a field, of degree two over $\QQ$. As $B\otimes_\QQ K\cong K^2$ we get $B\cong K$.
Note that $S_i\cong S_{m-i}$.

The interesting case is when $m$ is even and $2l=m$. 
Both summands of $S_{l,K}$ are in the same half spin representation (hence in the same eigenspace of the center of $C^+(V)_K$). The maps $R_f$ and $R_{\bar{f}}$ generate $End_{sl(V_+)}(S_{l,K})\cong M_2(K)$, hence 
$End_{su(H)}(S_l)=(\delta,d\delta)$ (see the proof of Prop.\ \ref{sodec}). 
If $d\equiv 0\,(4)$ this quaternion algebra
is a skew field and hence $S_l$ is irreducible. If $d\equiv 2\,(4)$, $S_l$ is irreducible in the non-split case and in the split case is isomorphic to 
$(S_l')^2$ with $S_{l,K}'\cong  \Lambda^l V_+$.
Thus we always have $End_{su(H)}(S_l)\cong D$, and $K\subset D$,
but $D$ is not
a skew field in the split case (then $D\cong M_2(\QQ)$). 
Since $S'_{l,K}$ is irreducible we have $End_{su(H)}(S'_l)\cong \QQ$.
\qed

\subsection{}
The previous propositions show that the Kuga-Satake Hodge structure
$(C^+(V),h_s)$ associated to $(V,h,\psi,K)$ decomposes as
$$
C^+(V)\cong (S_0\oplus S_1\oplus\ldots\oplus S_m)^{2^{m-2}}
$$
with $m=\dim_K V$ and they give the decomposition in simple factors
as well as the endomorphism rings
of the $S_i$ in the generic case. 

The Hodge structure on the $S_i$ can be obtained as follows.
Since $V$ is $K$-vector space, the exterior products $\wedge^i_K V$ are well-defined and, as $\QQ$-vector spaces, they have dimension $2{m\choose i}$,
which is just the dimension of the summand $S_i$. Weil already pointed out that
there is a natural inclusion 
$$
\wedge^i_K V\hookrightarrow \wedge^i V,
$$
and the $\wedge^i_K V$ are sub-Hodge structures of $\wedge^i V$ of weight $2i$.
Combining Tate and half twists of these, one obtains weight 1 Hodge structures which are the summands of the Kuga-Satake Hodge structure. Note that
the moduli of $V$, which has weight two, and of $S_1$ (with its $K$-action)
which has weight one, 
are both the $m-1$-ball.
Part of this theorem was already proved by Voisin in \cite{Vks}.

\subsection{Theorem}\label{ksdec}
{\it
 Let $(V,h,\psi)$ be a polarized weight $2$ Hodge structure with $\dim V^{2,0}=1$ which is of CM-type for an imaginary quadratic field $K$.

Then the Hodge structure on the summand $S_i$ (see Proposition \ref{sudec})
of the Kuga-Satake Hodge structure is:
$$
S_i\cong \left(\wedge_K^i V\right)(i-1)_\half.
$$
In particular, $S_1=V_\half$ and $S_0$ is the CM-type
Hodge structure of weight one on $K$.
Moreover:
$$
V\hookrightarrow S_0\otimes S_1.
$$
The dimensions of the eigenspaces for the $K$-action on $S_i^{1,0}$ are:
$$
\dim S_{i,\sigma}^{1,0}={m-1\choose i-1},\qquad 
\dim S_{i,\bar{\sigma}}^{1,0}={m-1\choose i},\qquad{\rm with}\quad
m=\dim_K V.
$$
}

\ts
Let $V_\CC=V_+\oplus V_-$ be the decomposition in eigenspaces for the $K$-action, combining this with the Hodge decomposition we get:
$$
V_\CC=V^{2,0}_+\oplus V^{1,1}_+\oplus V_-^{0,2}\oplus V_-^{1,1}.
$$  
We choose a basis $e_1,\ldots,e_m$ of $V_+$ and $e_{m+1},\ldots,e_n$ of $V_-$ such that 
$$
h(z)=diag(z^2,1\ldots,1,z^{-2},1,\ldots,1)\qquad(\in SO(\psi)(\CC))\qquad(z\in S^1),
$$
so, for example, $V^{0,2}=V^{0,2}_-=\langle e_{m+1}\rangle$.
Then $h(z)$ lies in the 1-parameter subgroup generated by 
$H_1:=\half(e_1\wedge e_{m+1})\in \wedge^2 V_\CC\cong so(\psi)_\CC$.
From the proof of \cite{FH}, 20.15, one finds that $H_1$ multiplies 
$w:=e_{i_1}\wedge\ldots\wedge e_{i_k}$, $1\leq i_1<\ldots<i_k\leq m$, by $+\half$
if $i_1=1$ and else by $-\half$. Hence $h_s(z)$ multiplies  $w$ by $z$
if $i_1=1$ and else by $\bar{z}$.
Since:
$$
S_{i,\CC}=\wedge^i V_+\oplus \wedge^{m-i} V_+=
\left(V^{2,0}_+ \otimes \wedge^{i-1}V_+^{1,1}\right)\oplus
\left(\wedge^i V_+^{1,1}\right)\oplus
\left(V^{2,0}_+ \otimes \wedge^{m-i-1}V_+^{1,1}\right)\oplus
\left(\wedge^{m-i} V_+^{1,1}\right)
$$
the action of $h_s(z)$ is by $diag(z,\bar{z},z,\bar{z})$.
The isomorphism $\wedge^{m-i}V_+\cong \wedge^iV_+^*\cong \wedge^iV_-$
induces
$$
V_+^{2,0}\otimes \wedge^{m-i-1} V_+^{1,1}\cong \wedge^i V_-^{1,1},\qquad
\wedge^{m-i}V_+^{1,1}\cong V_-^{2,0}\otimes \wedge^{i-1}V^{1,1}_-.
$$ 
Therefore the Hodge decomposition of $S_{i,\CC}\cong \wedge^iV_+\oplus
\wedge^i V_-$ is given by:
$$
S_i^{1,0}=\left(V^{2,0}_+ \otimes \wedge^{i-1}V_+^{1,1}\right)\oplus
\left(\wedge^{i}V_-^{1,1}\right).
$$
The eigenspace decomposition of $S^{1,0}_i$ is now obvious and we see that 
$S_1=V_\half$.

On the other hand, following Weil, we have
$$
(\wedge_K^i V)\otimes_\QQ\CC=
\left(\wedge^iV_+\right)\oplus \left(\wedge^iV_-\right)\quad
\hookrightarrow \quad   
\oplus_j\left( \wedge^{i-j}V_+\otimes\wedge^jV_-\right) 
= (\wedge^iV)\otimes\CC.
$$
The Hodge structure on $\wedge_K^iV$ induced from this inclusion is 
$$
(\wedge_K^i V)\otimes_\QQ\CC=
\left(V^{2,0}_+\otimes \wedge^{i-1}V^{1,1}_+\right)
\oplus
\left(\wedge^i V_+^{1,1}\oplus \wedge^iV_-^{1,1}\right) \oplus
\left(\wedge^{i-1}V^{1,1}_+\otimes V^{0,2}_-\right),
$$
thus the Hodge numbers are $(2,0)+(i-1,i-1)=(i+1,i-1)$, $(i,i)$ and $(i-1,i+1)$.
Therefore if we Tate twist $(i-1)$-times and then do a half twist we obtain a Hodge structure of weight one which is just the one obtained from $h_s$. 

The inclusion $V\subset V_\half\otimes K(-1)_\half=S_1\otimes S_0$
follows from Proposition \ref{halftens}.
\qed

\subsection{Example.}
We construct, geometrically,
a $9$ dimensional family of polarized Hodge structures with $h^{2,0}=1$, $h^{1,1}=18$
with CM by the field $K\cong\QQ(\sqrt{-3})$.
For a Hodge structure $V$ of this family
we identify the Hodge structures $S_0$ and $S_1$
as in Theorem \ref{ksdec} and we give a geometrical realisation
of the inclusion $V\hookrightarrow S_0\otimes S_1$.

For $a_1,\ldots,a_{12}\in \CC$ we define an (isotrivial) elliptic surface $S$ over $\PP^1$
by the the Weierstrass model:
$$
S:\quad Y^2=X^3+\prod_{i=1}^{12} (t-a_i),\qquad
S\rightarrow\PP^1,\quad(X,Y,t)\longmapsto t.
$$
Since $S$ has twelve fibers which are cuspidal it is a K3 surface
(and $\omega:=Y^{-1}{\rm d}X\wedge {\rm d}t$ is a no where zero holomorphic 2-form on $S$).

The orthogonal complement in $H^2(S,\QQ)$ of the classes of a fiber and 
the section at infinity is a sub-Hodge structure $V$ of dimension 20
in with $V^{2,0}=1$ and the field $K=\QQ(\sqrt{-3})$ acts on $V$
via the automorphism $(X,Y,t)\mapsto (\zeta^2 X,Y,t)$ with a primitive $6$-th root of unity $\zeta$.

Define curves $C$, $C'$, of genus $25$ and $1$ by:
$$
C:\qquad y^6=\prod_{i=1}^{12}(x-a_i);\qquad\qquad C':\qquad v^2=u^3-1.
$$
Both of these curves have automorphisms of order 6:
$$
\psi:C\rightarrow C,\quad (x,y)\mapsto (x,\zeta y);\qquad\qquad
\psi':C'\rightarrow C',\quad \psi':(u,v)\mapsto (\zeta^2 u, -v).
$$
The surface $S$ is the (minimal model of the desingularisation of the)
quotient of $C\times C'$ by the automorphism
 $\phi=(\psi^{-1},\psi')$ of order 6,
the quotient map is given by
$$
\pi:S=C\times C'\longrightarrow S,\qquad 
\Bigl((x,y),\,(u,v)\Bigr)\longmapsto (X,Y,t)=(y^2u,y^3v,x).
$$

To define $V_1$, consider the following rational 1-forms on $C$:
$$
\omega_{a,b}:=x^ay^b \frac{{\rm d} x}{y^5},\qquad{\rm note}\quad
\psi^*:\omega_{a,b}\mapsto \zeta^{b+1}\omega_{a,b}.
$$ 
It is easy to check that the $\omega_{a,b}$ with $a,\,b\geq 0$ and $a+2b\leq 8$
are a basis of $H^0(C,\omega_C)$. In particular, the eigenspace of $\psi^*$
with eigenvalue $\zeta$ has dimension $9$ (and is spanned by the $\omega_{a,0}$ with $0\leq a\leq 8$)
whereas the eigenspace with eigenvalue $\zeta^{-1}=\zeta^5$ has dimension 1
(and is spanned by $\omega_{0,4}$). 

Let $V_1\subset H^1(C,\QQ)$ be the $\QQ$-subspace
on which the eigenvalues of $\zeta$ are primitive $6$-th roots of unity.
Then $\dim V_1=20$, and the associated abelian variety is of Weil type $(1,9)$.
Let $V_0:=H^1(C',\QQ)$, note $\psi'$ acts 
on $H^0(C',\omega_{C'})=\langle {\omega':=\rm d}u/v\rangle$ as $\zeta^{-1}$.

The pull-back $\pi^*$ maps $V$ into the $\phi$-invariants in 
$H^1(C,\QQ)\otimes H^1(C',\QQ)$ and it is easy to verify
that these invariants are exactly
the $\phi$-invariants in $V_1\otimes V_0$.
For dimension reasons we then have:
$$
V\cong \pi^*V=\left( V_0\otimes_\QQ V_1\right)^{\langle \phi\rangle}.
$$
Since $V_0\cong K_{-\half}$, the half twist of this identity gives
$
V_\half\subset K\otimes V_1\cong V_1^{\oplus 2},
$
which implies that the half twist of $V$ is just $V_1$:
$$
V_\half\cong V_1
$$
and that $\pi^*$
is a geometrical realization of the Kuga-Satake correspondence.

The parameter space of 20 dimensional Hodge structures with CM by $K$
and $V^{2,0}=1$ is 
(a quotient of) the $9$-ball (cf.\ \ref{moduli}). 
The K3 surfaces in this example are parametrized by $12$ points in $\PP^1$,
Deligne and Mostov (\cite{DM}) actually showed that the geometrical quotient
$(\PP^1)^{12}//PGL(2)$ is a $9$-ball quotient.
See \cite{ravi} for old and new results on this moduli space.

\

\end{document}